\documentclass[12pt]{amsart}
\usepackage{amssymb}
\usepackage{amsfonts}
\usepackage{amsmath}
\usepackage{array}
\usepackage{rotating}
\usepackage{epsf}

\newtheorem{example}{Example}[section]

\newtheorem{note}[example]{Note}
\newtheorem{theorem}[example]{Theorem}

\newtheorem{corollary}[example]{Corollary}

\newtheorem{definition}[example]{Definition}

\def\Proof{\noindent \it Proof -- \rm}
\def\qed{\hspace{3.5mm} \hfill \vbox{\hrule height 3pt depth 2 pt width 2mm}
\bigskip}

\def\FQSym{{\bf FQSym}}
\def\WQSym{{\bf WQSym}}

\def\ssh{\Cup}

\def\Std{{\rm Std}}

\def\<{\langle}
\def\>{\rangle}

\def\F{{\bf F}}

\def\G{{\bf G}}

\def\M{{\bf M}}

\def\SG{{\mathfrak S}}

\def\goth{\mathfrak}
\def\Sym{{\bf Sym}}
\def\NCSF{{\bf Sym}}

\def\shuff#1#2{\mathbin{
\hbox{\vbox{ \hbox{\vrule \hskip#2 \vrule height#1 width 0pt
}%
\hrule}%
\vbox{ \hbox{\vrule \hskip#2 \vrule height#1 width 0pt
\vrule }%
\hrule}%
}}}

\def\binomial#1#2{\left(\,\begin{matrix}#1 \cr #2\end{matrix}\,\right)}
\def\binomial#1#2{\binom{#1}{#2}}

\def\gf#1#2{\genfrac{}{}{0pt}{}{#1}{#2}}

\def\shuf{{\mathchoice{\shuff{7pt}{3.5pt}}%
{\shuff{6pt}{3pt}}%
{\shuff{4pt}{2pt}}%
{\shuff{3pt}{1.5pt}}}}%
\def\shuffle{\,\shuf\,}

\def\Tabvrule{\vrule width-0.4pt}       % Difference de largeur
\def\Tabhrule{\hrule \hrule height-0.4pt} % Difference de hauteur
\def\Tabstrut{\vrule height2.2ex % Sur la ligne
                     depth0.8ex  % Sous la ligne
                     width0ex    % centrage horizontal
\relax}

\def\PasCase#1{\omit%
            $\vcenter{\hbox {\vbox to 0.4pt{}}
               \hbox{\makebox[3ex]{\Tabstrut$#1$}}}%
               \Tabvrule$}
\def\PasCasePoint{\PasCase{\cdot}}
\def\DessinCarre#1{%
    \vcenter{\hbox{}\hrule
             \hbox{\vrule\makebox[3ex]{\Tabstrut$#1$}\vrule}\Tabhrule}%
             \Tabvrule}
\def\GenRuban#1{\vcenter{\halign{&$\DessinCarre{##}$\cr#1}}\egroup}

\def\sTabvrule{\vrule width-0.4pt}
\def\sTabhrule{\hrule \hrule height-0.4pt}
\def\sTabstrut{\vrule height1.6ex depth0.6ex width0ex \relax}
\def\sDessinCarre#1{%
    \vcenter{\hbox{}\hrule
             \hbox{\vrule\makebox[2.3ex]%
                  {\sTabstrut$\scriptstyle#1$}\vrule}\sTabhrule}%
             \sTabvrule}
\def\sGenRuban#1{\vcenter{\halign{&$\sDessinCarre{##}$\cr#1}}\egroup}

\def\ruban{%
  \bgroup
  \let\ =\omit
  \let\\=\cr
  \let\x=\times
  \let\.=\PasCasePoint
  \offinterlineskip
  \GenRuban}

\def\sruban{%
  \bgroup
  \let\ =\omit
  \let\x=\times
  \let\\=\cr
  \offinterlineskip
  \sGenRuban}

\def\DesS{\operatorname{Des}}
\def\CDes{\operatorname{D}}
\def\GDes{\operatorname{GDes}}
\def\PW{{\rm PW}}   % packed words = mots tasses
\def\K{{\mathbb K}}   % corps de base
\def\convW{{*_W}}       % une loi de convolution
\def\pack{\operatorname{pack}}
\def\NCQSym{{\bf NCQSym}}      % Mots initiaux (pour Bergeron-Zabrocki)
\def\GC{\operatorname{GC}}
\def\WC{\operatorname{WC}}
\def\TE{{\bf T}}
\def\TEb{{\bf T'}}
\def\cal{\mathcal}

\title[Permutations, monomial NCSF and Genocchi numbers]%
{Permutation statistics related to a class\\
 of noncommutative symmetric functions\\
and generalizations of the Genocchi numbers}

\author[F.~Hivert, J.-C.~Novelli, L.~Tevlin, J.-Y.~Thibon]%
{Florent Hivert, Jean-Christophe Novelli,\\ Lenny Tevlin, and Jean-Yves Thibon}

\address[F. Hivert]{LITIS, Universit\'e de Rouen ; Avenue de l'universit\'e ;
76801 Saint \'Etienne du Rouvray, France\\}
\address[J.-C Novelli, J.-Y. Thibon] {Institut Gaspard Monge, Universit\'e
Paris-Est Marne-la-Vall\'ee \\
5 Boulevard Descartes \\Champs-sur-Marne \\77454 Marne-la-Vall\'ee cedex 2 \\
France}
\address[L. Tevlin]{Physics Department, Yeshiva University, 
                    500 West 185th Street, New York, N.Y. 10033, USA}

\email[F. Hivert]{hivert@univ-rouen.fr}
\email[J.-C. Novelli]{novelli@univ-mlv.fr}
\email[L. Tevlin]{tevlin@yu.edu}
\email[J.-Y. Thibon]{jyt@univ-mlv.fr} 
\date{\today}

\begin{document}

\begin{abstract}
We prove conjectures of the third author~[L. Tevlin, Proc. FPSAC'07, Tianjin]
on two new bases of noncommutative symmetric functions: the transition
matrices from the ribbon basis have nonnegative integral coefficients.  This
is done by means of two composition-valued statistics on permutations and
packed words, which generalize the combinatorics of Genocchi numbers.
\end{abstract}

\keywords{Noncommutative symmetric functions, quasideterminants, permutation
statistics, Genocchi numbers}

\subjclass[2000]{05E05;15A15; 16W30}

\maketitle

%%%%%%%%%%%%%%%%%%%%%%%%%%%%%%%%%%%%%%%%%%%%%%%%%%%%%%%%%%%%%%%%%%%%%%%%%%%%%%%
%%%%%%%%%%%%%%%%%%%%%%%%%%%%%%%%%%%%%%%%%%%%%%%%%%%%%%%%%%%%%%%%%%%%%%%%%%%%%%%
%%%%%%%%%%%%%%%%%%%%%%%%%%%%%%%%%%%%%%%%%%%%%%%%%%%%%%%%%%%%%%%%%%%%%%%%%%%%%%%

\section{Introduction}

In the theory of noncommutative symmetric functions \cite{NCSF1}, the
self dual commutative Hopf algebra $Sym$ of ordinary symmetric functions
is replaced by a pair of mutually dual Hopf algebras
$(\Sym, QSym)$, respectively called Noncommutative Symmetric Functions,
and Quasi-symmetric functions. The usual bases of $Sym$ are usually lifted
only on one side (with the notable exception of Schur functions, which
admit natural analogs on both sides). In particular, monomial symmetric
functions $m_\mu$ split into the quasi-monomial functions $M_I$ on the
quasi-symmetric side, and their dual basis $h_\mu$ is lifted on the
noncommutative side, in the form of the homogeneous products $S^I$.

In \cite{Tev}, the third author has proposed a construction of noncommutative
monomial and forgotten symmetric functions, and conjectured positivity
properties of certain transition matrices involving the new bases. The purpose
of the present article is to prove these conjectures, by providing
combinatorial interpretations.

These interpretations rely on new permutations statistics, which generalize
the combinatorics related to the Genocchi numbers.

{\footnotesize \it Notations. \rm We shall depart from the notation of
\cite{Tev} and write $\Psi_I$ instead of $M^I$, and $L_I$ instead of $L^I$.
Other notations are as in~\cite{NCSF1}.
See \cite{GGRW} for background on quasideterminants.
}

{\footnotesize {\it Acknowledgments.-}
This work has been partially supported by Agence Nationale de la Recherche,
grant ANR-06-BLAN-0380.
The authors would also like to thank the contributors of the MuPAD project,
and especially those of the combinat package, for providing the development
environment for this research (see~\cite{HT} for an introduction to
MuPAD-Combinat).
}
%%%%%%%%%%%%%%%%%%%%%%%%%%%%%%%%%%%%%%%%%%%%%%%%%%%%%%%%%%%%%%%%%%%%%%%%%%%%%%%
%%%%%%%%%%%%%%%%%%%%%%%%%%%%%%%%%%%%%%%%%%%%%%%%%%%%%%%%%%%%%%%%%%%%%%%%%%%%%%%
%%%%%%%%%%%%%%%%%%%%%%%%%%%%%%%%%%%%%%%%%%%%%%%%%%%%%%%%%%%%%%%%%%%%%%%%%%%%%%%
\section{Background}

%%%%%%%%%%%%%%%%%%%%%%%%%%%%%%%%%%%%%%%%%%%%%%%%%%%%%%%%%%%%%%%%%%%%%%%%%%%%%%%
\subsection{Noncommutative symmetric functions}

Recall that the algebra $\Sym$ of noncommutative symmetric functions
is a graded free associative algebra over a sequence $S_n$ of indeterminates,
with ${\rm deg}\, S_n=n$. Among other sequences of generators,
the {\em noncommutative power sums of the first kind} $\Psi_n$ are
defined by an oriented analog of Newton's recursion, which may be
solved in terms of quasideterminants \cite{NCSF1,GGRW}.
The following definition \cite{Tev} refines formulas (39) and (40)
of \cite{NCSF1}, and defines an analog of the monomial basis which
extends the $\Psi_n$.

\begin{definition}
The \textbf{noncommutative monomial symmetric function} corresponding to a
composition $ I = (i_1, \dots, i_r)$ is defined as a quasideterminant of an
$r$ by $r$ matrix:
\label{def-primm}
\begin{equation}
 r \Psi_I \equiv r \Psi_{(i_1, \ldots, i_r)} = (-1)^{r - 1}
\begin{vmatrix}
\Psi_{i_r} & 1 & 0 & \dots & 0& 0\\
\Psi_{ i_{n -1} + i_r} & \Psi_{i_{n - 1}} & 2 & \ldots & 0 & 0 \\
\vdots &  \vdots & \vdots & \vdots &  \vdots & \vdots \\
\Psi_{i_2 + \ldots + i_r} & \ldots & \ldots & \ldots&  \Psi_{i_2}& n - 1 \\
\fbox{$\Psi_{i_1 + \ldots + i_r}$} & \ldots & \ldots & \ldots &\Psi_{i_1 +
i_2}& \Psi_{i_1}
\end{vmatrix}
\end{equation}
\end{definition}
where $r$ is the length of $ I $.
In particular,
\begin{equation}
\Psi_{(n)}=\Psi_n,
\text{\ and\ }
\Psi_{1^r} = \Lambda_r
\end{equation}
where $ \Lambda_r $ is an elementary symmetric function.

The quasideterminants may be recursively evaluated by means of the
following generalized Newton relations:

\begin{equation}
\label{eq-linear}
\begin{split}
 r \Psi_{i_1, \dots, i_r} &= \Psi_{i_1}\Psi_{i_2, \dots , i_r} - \Psi_{i_1 +
i_2}\Psi_{i_3, \dots, i_r}
+ \dots\\
& + (-1)^{s - 1} \Psi_{i_1 + \dots + i_s} \Psi_{i_{s +1},  \dots, i_r} +
\dots +
(-1)^r \Psi_{i_1 + \dots + i_r} .
\end{split}
\end{equation}

From a noncommutative analog of the quasi-monomial basis $M_I$,
one can define an analog of Gessel's fundamental basis $F_I$ by
\begin{equation}
\label{LPsi}
L_I=\sum_{J\succeq I}\Psi_J.
\end{equation}
Define the coefficients $G_{IJ}$ by the expansion
\begin{equation}
R_I=\sum_J G_{IJ} L_J.
\end{equation}
It has been conjectured in \cite{Tev} that these numbers
are nonnegative integers. Our aim is to prove this fact
by means of a combinatorial interpretation.

%%%%%%%%%%%%%%%%%%%%%%%%%%%%%%%%%%%%%%%%%%%%%%%%%%%%%%%%%%%%%%%%%%%%%%%%%%%%%%%
\subsection{Free quasi-symmetric functions}

Let us fix an infinite ordered alphabet $A=\{a_1<\dots<a_n<\dots\}$.
The \emph{standardized word} $\Std(w)$ of a word $w\in A^*$ is the permutation
obtained by iteratively scanning $w$ from left to right, and labelling
$1,2,\ldots$ the occurrences of its smallest letter, then numbering the
occurrences of the next one, and so on.

With a permutation $\sigma$, we associate the polynomial
\begin{equation}
\F_\sigma :=\sum_{\Std(w)=\sigma^{-1}}w\,.
\end{equation}
These polynomials span a subalgebra of $\K\langle A\rangle$, called $\FQSym$
for Free Quasi-Sym\-metric functions~\cite{NCSF6}. Note that the field $\K$ is
assumed to be of characteristic zero.
Their product rule is given by
\begin{equation} 
\label{prodF-fq}
\F_{\sigma'} \F_{\sigma''} = \sum_{\sigma \in \sigma'\ssh \sigma''}
\F_\sigma\,,
\end{equation}
where the \emph{shifted shuffle} $\sigma'\ssh \sigma''$ of two packed words
is defined as
\begin{equation} 
\sigma'\ssh \sigma'' = \sigma'\shuffle\sigma''[|\sigma|],
\end{equation}
the $k$-shift $w[k]$ of a word $w$ being obtained by replacing each letter
$w_i$ by $w_i+k$, and $\shuffle$ is the usual shuffle product on words defined
recursively by
\begin{equation}
(au)\shuffle (bv) = a\cdot(u\shuffle (bv)) + b\cdot ((au)\shuffle v),
\end{equation}
with $u\shuffle\epsilon=\epsilon\shuffle u=u$ if $\epsilon$ is the empty word.

We shall make use of the basis $\G_\sigma$ of $\FQSym$, dual to $\F_\sigma$,
defined by $\G_\sigma:=\F_{\sigma^{-1}}$.

%%%%%%%%%%%%%%%%%%%%%%%%%%%%%%%%%%%%%%%%%%%%%%%%%%%%%%%%%%%%%%%%%%%%%%%%%%%%%%%
\subsection{Word quasi-symmetric functions}

The \emph{packed word} $u=\pack(w)$ associated with a word $w\in A^*$ is
obtained by the following process. If $b_1<b_2<\ldots <b_r$ are the letters
occuring in $w$, $u$ is the image of $w$ by the homomorphism
$b_i\mapsto a_i$.
A word $u$ is said to be \emph{packed} if $\pack(u)=u$. We denote by $\PW$ the
set of packed words.
With such a word, we associate the polynomial
\begin{equation}
\M_u :=\sum_{\pack(w)=u}w\,.
\end{equation}
These polynomials span a subalgebra of $\K\langle A\rangle$, called $\WQSym$
for Word Quasi-Symmetric functions~\cite{Hiv,NT06} (and called $\NCQSym$
in~\cite{BZ}), the invariants of the noncommutative quasi-symmetrizing action.
Their product rule is given by
\begin{equation} 
\label{prodG-wq}
\M_{u'} \M_{u''} = \sum_{u \in u'\convW u''} \M_u\,,
\end{equation}
where the \emph{convolution} $u'\convW u''$ of two packed words
is defined as
\begin{equation} 
u'\convW u'' =
\sum_{\genfrac{}{}{0pt}{}{v,w ; u=v\cdot w\,\in\,\PW}
{ \pack(v)=u', \pack(w)=u''}} u\,.
\end{equation}

%%%%%%%%%%%%%%%%%%%%%%%%%%%%%%%%%%%%%%%%%%%%%%%%%%%%%%%%%%%%%%%%%%%%%%%%%%%%%%%
%%%%%%%%%%%%%%%%%%%%%%%%%%%%%%%%%%%%%%%%%%%%%%%%%%%%%%%%%%%%%%%%%%%%%%%%%%%%%%%
%%%%%%%%%%%%%%%%%%%%%%%%%%%%%%%%%%%%%%%%%%%%%%%%%%%%%%%%%%%%%%%%%%%%%%%%%%%%%%%
\section{A statistic on permutations generalizing Genocchi numbers}

Genocchi numbers (sequence A001469 of~\cite{Slo}) are known to count a large
variety of combinatorial objects, among which numerous sets of permutations.
Our statistic derives directly from the most classical of those sets: it is
the number of permutations of $\SG_{2n}$ such that each even integer is
followed by a smaller integer and each odd integer is either followed by a
greater one, or at the last position of the permutation.

Let us define the \emph{Genocchi descent set} (G-descent set for short) of a
permutation $\sigma\in\SG_n$ as
\begin{equation}
\GDes(\sigma) := \{i\in[2,n] | \sigma_{j}=i
                 \Longrightarrow \sigma_{j+1}<\sigma_j
\}.
\end{equation}
In other words, $\GDes(\sigma)$ is the set of \emph{values} of the descents of
$\sigma$, different from the usual set $\DesS(\sigma)$ which records the
\emph{positions} of the descents of $\sigma$. Astonishingly enough,
this G-statistic behaves very differently from the classical descent
statistic.
From the G-descent set, we define the \emph{Genocchi composition of descents}
(or G-compos\-ition, for short) $\GC(\sigma)$ of a permutation, as the integer
composition $I$ of $n$ whose descent set is $\{d-1|d\in\GDes(\sigma)\}$.

The following tables represent the G-composition of all permutations of
$\SG_2$, $\SG_3$, and $\SG_4$.

\begin{equation}
\begin{array}{|c|c|}
\hline
2 & 11 \\
\hline
\hline
12 & 21 \\
\hline
\end{array}
\hskip2cm
\begin{array}{|c|c|c|c|}
\hline
3   & 21  & 12  & 111\\
\hline
\hline
123 & 132 & 213 & 321\\
    & 231 &     &    \\
    & 312 &     &    \\
\hline
\end{array}
\end{equation}

\begin{equation}
\begin{array}{|c|c|c|c|c|c|c|c|}
\hline
4    & 31   & 22   & 211  & 13   & 121  & 112  & 1111\\
\hline
\hline
1234 & 1243 & 1324 & 1432 & 2134 & 2143 & 3214 & 4321 \\
     & 1342 & 2314 & 2431 &      & 3421 &      &      \\
     & 1423 & 3124 & 3142 &      & 4213 &      &      \\
     & 2341 &      & 3241 &      &      &      &      \\
     & 2413 &      & 4132 &      &      &      &      \\
     & 3412 &      & 4231 &      &      &      &      \\
     & 4123 &      & 4312 &      &      &      &      \\
\hline
\end{array}
\end{equation}
More combinatorial properties of these numbers, including a hook-length
formula will be given in~\cite{NTW}.

%%%%%%%%%%%%%%%%%%%%%%%%%%%%%%%%%%%%%%%%%%%%%%%%%%%%%%%%%%%%%%%%%%%%%%%%%%%%%%%
%%%%%%%%%%%%%%%%%%%%%%%%%%%%%%%%%%%%%%%%%%%%%%%%%%%%%%%%%%%%%%%%%%%%%%%%%%%%%%%
%%%%%%%%%%%%%%%%%%%%%%%%%%%%%%%%%%%%%%%%%%%%%%%%%%%%%%%%%%%%%%%%%%%%%%%%%%%%%%%
\section{A $\NCSF$ quotient of $\FQSym$}

Let $\sim$ be the equivalence relation defined by $\sigma\sim\tau$
iff $\GC(\sigma)=\GC(\tau)$.
Let ${\cal J}$ be the subspace of $\FQSym$ spanned by the differences
%Let us consider the algebra $\TE$ as the quotient of $\FQSym$ by the relations
%\F_\sigma \equiv \F_\tau
%\Longleftrightarrow
%\GC(\sigma)=\GC(\tau).
\begin{equation}
\{ \F_\sigma - \F_\tau | \sigma\sim\tau\}.
\end{equation}

\begin{theorem}
\label{TEQuot}
${\cal J}$ is a two-sided ideal of $\FQSym$, and the quotient
$\TE=\FQSym/{\cal J}$ is isomorphic to $\NCSF$ as an algebra.

Moreover, let $T_I$ be the image in $\TE$ of the $\F_\sigma$ such that
$\GC(\sigma)=I$.
Then
\begin{equation}
\label{prodT}
T_I T_J =
\sum_{K} C_{I,J}^K T_K,
\end{equation}
where $C_{I,J}^K$ is computed as follows.
Let $K'$ and $K''$ be the compositions such that $|K'|=|I|$ and either
$K=K'\cdot K''$, or $K=K'\triangleright K''$.
If $K'$ is not coarser than $I$ or if $K''$ is not finer than $J$, then
$C_{I,J}^K$ is $0$.
Otherwise,
\begin{equation}
\label{cijk}
C_{I,J}^K= \binomial{|I|+l(J)-l(I)}{l(K)-l(I)}
\end{equation}
\end{theorem}

\Proof
We have to prove that the set (with multiplicities) of G-compositions of the
shifted shuffle of two permutations depends only on the G-compositions
of the permutations.

Let $\sigma\in\SG_m$ and $\tau\in\SG_{n}$ and let $I$ and $J$ be their
respective G-compositions.
Let $K=(k_1,\dots,k_r)$ be a composition of $m+n$ and let us compute the
number of permutations $\mu$ in $\sigma\ssh\tau$ such that $\GC(\mu)=K$.
We shall need the unique compositions $K'$ and $K''$ such that
$K=K'\cdot K''$ or $K=K'\triangleright K''$, with $|K'|=|I|$.

Let us now consider which letters can follow the letters $x$ from 1 to $m$ in
$\mu$.
We have four cases:
\begin{itemize}
\item[(1a)] $x$ is a G-descent of $\mu$ and is not a G-descent of $\sigma$,
\item[(2a)] $x$ is a G-descent of $\mu$ and is a G-descent of $\sigma$,
\item[(3a)] $x$ is not a G-descent of $\mu$ and is not a G-descent of $\sigma$,
\item[(4a)] $x$ is not a G-descent of $\mu$ and is a G-descent of $\sigma$.
\end{itemize}

The first case implies that $K$ cannot be the G-composition of a word in the
shifted shuffle of $\sigma$ and $\mu$. Let us now restrict to compositions $K$
such that $K'$ is coarser than $I$.
The second case implies that $x$ has to be followed in $\mu$ by a letter
coming from $\sigma$.
The third one implies nothing about $x$.
The fourth one implies that $x$ has to be followed in $\mu$ by a letter coming
from $\tau$.

Let $f(\sigma)$ be the number of occurrences of the fourth case.
Let $g(\sigma)$ be the number of occurrences of the third and fourth cases,
plus one.

Let us now consider which letters can follow the letters $x$ from $m+1$ to
$m+n$ in $\mu$.
We have again four cases:
\begin{itemize}
\item[(1b)] $x$ is a G-descent of $\mu$ and is not a G-descent of $\tau[m]$,
\item[(2b)] $x$ is a G-descent of $\mu$ and is a G-descent of $\tau[m]$,
\item[(3b)] $x$ is not a G-descent of $\mu$ and is not a G-descent of $\tau[m]$,
\item[(4b)] $x$ is not a G-descent of $\mu$ and is a G-descent of $\tau[m]$.
\end{itemize}

The first case implies that $x$ has to be followed in $\mu$ by a letter coming
from $\sigma$. 
The second case implies nothing about $x$.
The third one implies that $x$ has to be followed in $\mu$ by a letter coming
from $\tau$.
The fourth one implies that $K$ cannot be the G-composition of a word in the
shifted shuffle of $\sigma$ and $\mu$. We now restrict to compositions $K$
such that $K''$ is finer than $J$.

This preliminary analysis proves that the number of permutations $\mu$ with
G-composition $K$ is equal to the number of ways of separating the letters of
$\tau$ in any number of blocks with given necessary separations (case 1b)
and necessary non-separations (case 3b) and put those blocks in the middle of
blocks of letters of $\sigma$, themselves separated into this number of blocks
with given necessary separations (case 4a), and necessary non-separations
(case 2a).
The \emph{number} of such blocks for each pair of permutations depends only on
the lengths of their G-compositions, and \emph{a fortiori} only on their
G-compositions, so that our equivalence relation on permutations indeed
induces a quotient algebra of $\FQSym$.

Let us now determine the structure constants of this algebra.
The previous remark can be reformulated as follows.
We have two cases depending on whether $K=K'\cdot K''$ or
$K=K'\triangleright K''$.

In the first case, the number of permutations $\tau$ with G-composition $I$ is
\begin{equation}
\sum_{k=1}^{\max(n,m)}
 \binomial{m-l(I)}{k-1-(l(I)-l(K'))}
 \binomial{n}{l(K'')+1-k}.
\end{equation}
This sum of binomial coefficients is easily simplified, and one gets
\begin{equation}
\binomial{m-l(I)+l(J)}{l(K')+l(K'')-l(I)}.
\end{equation}
In the second case, the number of permutations $\tau$ with G-composition $I$ is
\begin{equation}
\sum_{k=1}^{\max(n,m)}
 \binomial{m-l(I)}{k-1-(l(I)-l(K')}
 \binomial{n}{l(K'')-k}.
\end{equation}
Similarly, this sum of binomial coefficients reduces to
\begin{equation}
\binomial{m-l(I)+l(J)}{l(K')+l(K'')-1-l(I)},
\end{equation}
so that $C_{IJ}^K$ is indeed given by~(\ref{cijk}) in any case.

%This ends the computation of the $C_{I,J}^K$.
One can notice that these coefficients coincide, in the special case $I=(n)$
with those of the product $L_I L_J$ (see Proposition 4.6 of~\cite{Tev}), so
that, since the $L_n$ are algebraic generators of $\NCSF$, the $T_n$ are
algebraic generators of $\TE$.
Moreover, since $\NCSF$ is free over the sequence $L_n=S_n$, the algebra $\TE$
is free over the $T_n$.
\qed

\begin{example}
{\rm
Let $I=(2,2,1)$, $J=(1,3)$, and $K=(4,2,1,1,1)$.
We can choose $\sigma=32514$ and $\tau=2134$.

We then have $K'=(4,1)$ and $K''=(1,1,1,1)$. The coefficient of $K$ in
$T_I T_J$ is $\binom{5+2-3}{5-3}=6$ and, indeed, there are six permutations
in the shifted shuffle $32514\ssh 2134$ with G-composition $K$:
\begin{equation}
372685194,\ 376825194,\
376829514,\ 736825194,\
736829514,\ 768392514.
\end{equation}
Those six permutations are obtained as follows:
$\sigma$ has one necessary separation between $3$ and $2$ and one necessary
non-separation between $1$ and $5$, and nothing after~$4$.
The permutation $\tau[5]$ has two necessary separations, between $8$ and $9$,
and after $9$, and one necessary non-separation between $6$ and $8$. 
Then one inserts the blocks of $\tau[5]$ in $\sigma$, satisfying the
separation/non-separation constraints and gets the six permutations.
}
\end{example}

\begin{note}
{\rm
Note that this quotient is not a Hopf quotient, since as one can easily check,
${\cal J}$ is not a coideal. For example, $231\sim 312$ but
\begin{equation}
\overline\Delta(\F_{231}) = \F_{12}\otimes \F_{1} + \F_{1}\otimes\F_{21},
\text{\ and\ }
\overline\Delta(\F_{312}) = \F_{21}\otimes \F_{1} + \F_{1}\otimes\F_{12}.
\end{equation}
%It is not also a sub-coalgebra of $\FQSym$ on the $\G$ basis since
%\begin{equation}
%%\overline\Delta(\G_{231}) = \G_{21}\otimes \G_{1} + \G_{1}\otimes\G_{12},
%\overline\Delta(\G_{231}) = \G_{21}\otimes \G_{1} + \G_{1}\otimes\G_{12},
%\text{\ and\ }
%\overline\Delta(\G_{312}) = \G_{12}\otimes \G_{1} + \G_{1}\otimes\G_{21},
%\end{equation}
}
\end{note}

%%%%%%%%%%%%%%%%%%%%%%%%%%%%%%%%%%%%%%%%%%%%%%%%%%%%%%%%%%%%%%%%%%%%%%%%%%%%%%%
%%%%%%%%%%%%%%%%%%%%%%%%%%%%%%%%%%%%%%%%%%%%%%%%%%%%%%%%%%%%%%%%%%%%%%%%%%%%%%%
%%%%%%%%%%%%%%%%%%%%%%%%%%%%%%%%%%%%%%%%%%%%%%%%%%%%%%%%%%%%%%%%%%%%%%%%%%%%%%%
\section{Change of bases in $\NCSF$}

Thanks to the previous result, we have a map going from $\NCSF$ to itself in
a very unusual way: start with the injection of $\NCSF$ into $\FQSym^*$, and
compose it with the self duality isomorphism of $\FQSym$, which reads
\begin{equation} 
R_I := \sum_{\CDes(\sigma)=I} \G_\sigma = \sum_{\CDes(\tau^{-1})=I} \F_\tau,
\end{equation}
where $\CDes$ is the composition whose descent set is equal to the descent set
of $\sigma$, and then go from $\FQSym$ to $\NCSF$ by the G-quotient
homomorphism.

Let $\phi$ be the composition of those maps and let $R'_I$ be the image of
$R_I$ by $\phi$:
\begin{equation}
\begin{array}{cccc}
\phi:& \NCSF & \to     & \TE  \\
     & R_I   & \mapsto & R'_I.\\
\end{array}
\end{equation}
By definition of $\phi$, we have
\begin{equation}
\label{RpT}
R'_I := \sum_{\CDes(\sigma^{-1})=I} \overline{\F_\sigma}
=
\sum_{\genfrac{}{}{0pt}{}%
{\CDes(\sigma^{-1})=I}{\GC(\sigma)=J}} T_J,
\end{equation}
where  $\overline\F_\sigma$ is the image of $\F_\sigma$ by the G-quotient
homomorphism.
%so that
%\begin{equation}
%R'_I = \sum_{J\vDash n} G'_{IJ} T_J,
%\end{equation}
%
Then, since $L_n=R_n$ and $R'_n=\overline{\F_{12\dots n}}=T_n$, we have
$\phi(L_n)=T_n$ for all $n$, so that, thanks to the product formulas of $L_n$
and $T_n$, $\phi(L_I)=T_I$ for all compositions $I$.

Since the $T_n$ are algebraic generators of $\TE$, the algebra morphism $\phi$
is an isomorphism of algebras, so that, applying $\phi^{-1}$ to
Equation~(\ref{RpT}), one gets

\begin{theorem}
\label{ThRL}
Let $I$ be a composition of $n$. Then
\begin{equation}
\label{R2L}
R_I = \sum_{J\vDash n} G_{IJ} L_J,
\end{equation}
where $G_{IJ}$ is the number of permutations satisfying
$\CDes(\sigma^{-1})=I$ and $\GC(\sigma)=J$.
In particular, the $G_{IJ}$ are nonnegative integers.
\end{theorem}

Examples of the transition matrices are given in Section~\ref{tabsg},
together with the same matrices filled with the corresponding permutations.
The Genocchi numbers appear as the sums of the values in the rows indexed by
compositions of the form $(2^n)$ or $(2^n1)$.

Combining this last result with Equation~(\ref{LPsi}), one then gets

\begin{corollary}
\label{CorRPsi}
Let $I$ be a composition of $n$. Then
\begin{equation}
R_I = \sum_{J\vDash n} K_{IJ} \Psi_J,
\end{equation}
where $K_{IJ}$ are nonnegative integers. 
\end{corollary}

One can easily describe those integers in terms of permutations.
They can be described in a much more natural way in terms of packed words as
one shall see in the following section.

\begin{note}
{\rm
Theorems~\ref{prodT} and~\ref{ThRL} prove Conjecture 4.1 of~\cite{Tev}.
Theorem~\ref{ThRL} proves Conjecture 5.4 of~\cite{Tev}.
Corollary~\ref{CorRPsi} proves Conjecture 5.3 of~\cite{Tev}.
Corresponding statements for the ``forgotten'' basis are obtained by applying
the canonical involution $\omega$.
}
\end{note}

%\begin{equation}
%\end{equation}

%%%%%%%%%%%%%%%%%%%%%%%%%%%%%%%%%%%%%%%%%%%%%%%%%%%%%%%%%%%%%%%%%%%%%%%%%%%%%%%
%%%%%%%%%%%%%%%%%%%%%%%%%%%%%%%%%%%%%%%%%%%%%%%%%%%%%%%%%%%%%%%%%%%%%%%%%%%%%%%
%%%%%%%%%%%%%%%%%%%%%%%%%%%%%%%%%%%%%%%%%%%%%%%%%%%%%%%%%%%%%%%%%%%%%%%%%%%%%%%
\section{Permutations replaced by packed words}

The previous section was devoted to the study of the transition matrices from
$R$ to $L$, we now apply a similar analysis to the transition matrices from
$R$ to $\Psi$.
As already mentioned, the latter can be described in terms of the former,
since there is a very simple transition matrix from $L$ to $\Psi$: it is the
matrix of the refinement order on compositions.
Nevertheless, as the sum of the entries of the transition matrix $M_n(R,L)$ is
$n!$, the sum of the entries of a transition matrix $M_n(R,\Psi)$ is the
$n$th ordered Bell number (sequence A000670 of~\cite{Slo}) counting,
for example, set compositions (ordered set partitions), or packed words.

This suggests the existence of two statistics on packed words giving back
the entries of the transition matrices, exactly as in the $(R,L)$ case.
The algebraic context is essentially the same as before if one replaces the
algebra $\FQSym$ by the algebra $\WQSym$ (see~\cite{Hiv,NT06}).

The proof of the connection between the two statistics and the matrices
$M(R,\Psi)$ follows the same guidelines as the previous proof. We first define
a composition-valued statistic on packed words, then prove that this statistic
defines a quotient of $\WQSym$, isomorphic to $\NCSF$ as an algebra. Then,
comparing the structure constants of the natural base with those of the
$\Psi_I$, we prove that they are mapped to each other by a simple isomorphism,
hence giving the coefficients of the matrix $M(R,\Psi)$.

%%%%%%%%%%%%%%%%%%%%%%%%%%%%%%%%%%%%%%%%%%%%%%%%%%%%%%%%%%%%%%%%%%%%%%%%%%%%%%%
\subsection{A statistic on packed words}

Let $w$ be a packed word.
The \emph{Word composition} (W-comp\-os\-ition) of $w$ is the composition
whose descent set is given by the positions of the last occurrences of each
letter in $w$.
For example,
\begin{equation}
\WC(1543421323) = (2,3,2,2,1).
\end{equation}
Indeed, the descent set is $\{2,5,7,9,10\}$ since
the last $5$ is in position $2$,
the last $4$ is in position $5$,
the last $1$ is in position $7$,
the last $3$ is in position $9$,
the last $2$ is in position $10$.

The following tables represent the W-compositions of all packed words in
$\PW_2$ and $\PW_3$. One can recover from the matrix ${\goth M'}_4$ in
Section~\ref{tabsk} the W-compositions in $\PW_4$: it is the composition
indexing their row.

\begin{equation}
\begin{array}{|c|c|}
\hline
2 & 11 \\
\hline
\hline
11 & 12 \\
   & 21 \\
\hline
\end{array}
\hskip2cm
\begin{array}{|c|c|c|c|}
\hline
3   & 21  & 12  & 111\\
\hline
\hline
111 & 112 & 122 & 123\\
    & 121 & 211 & 132\\
    & 212 &     & 213\\
    & 221 &     & 231\\
    &     &     & 312\\
    &     &     & 321\\
\hline
\end{array}
\end{equation}

%%%%%%%%%%%%%%%%%%%%%%%%%%%%%%%%%%%%%%%%%%%%%%%%%%%%%%%%%%%%%%%%%%%%%%%%%%%%%%%
\subsection{A $\NCSF$ quotient of $\WQSym$}

Let $\sim$ be the equivalence relation on packed words defined by
$u\sim v$ iff $\WC(u)=\WC(v)$.
Let ${\cal J'}$ be the subspace of $\WQSym$ spanned by the differences
\begin{equation}
\{ \M_u - \M_v\, |\, u\sim v\}.
\end{equation}

\begin{theorem}
\label{TEQuotb}
${\cal J'}$ is a two-sided ideal of $\WQSym$, and the quotient
$\TEb$ defined by $\TEb=\WQSym/{\cal J'}$ is isomorphic to $\NCSF$ as an
algebra.

Moreover, let $U_I$ be the image of $\M_u$ in $\TEb$.
Then
\begin{equation}
\label{prodU}
U_I U_J :=
\sum_{K} D_{I,J}^K U_K,
\end{equation} 
where $D_{I,J}^K$ is computed as follows.
Let $K'$ and $K''$ be the compositions such that $|K'|=|I|$ and either
$K=K'\cdot K''$, or $K=K'\triangleright K''$.
If $K'$ is not coarser than $I$, then $D_{I,J}^K$ is $0$.
Otherwise,
\begin{equation}
\label{dijk}
D_{I,J}^K= \binomial{l(K)}{l(I)}.
\end{equation}
\end{theorem}

\Proof
The proof follows essentially the same lines as the proof of
Theorem~\ref{TEQuot}, so we only sketch it. In fact, the details
are much simpler than for Theorem~\ref{TEQuot}.
Looking at the definitions of the W-composition and of the convolution of
packed words, it is clear that the multiset of the W-compositions of the words
in the convolution of two packed words depends only on the W-compositions
of the words.
% Indeed, the convolution only amounts to choose which
%letters have to be equal between the left and the right word. So we can 
% So we obtain a quotient of $\WQSym$.
So the product is well-defined and $\TEb$ is a quotient of $\WQSym$.

Let us now see why the product $U_I U_J$ is given by Equations~(\ref{prodU})
and~(\ref{dijk}).
Let us choose two words $u$ and $v$ such that $\WC(u)=I$, and $\WC(v)=J$.
Since there is exactly one nondecreasing word having a given $\WC$, we
can assume that $u$ and $v$ are nondecreasing.
Let $|v|$ be the size of $v$.
Let us compute $\M_u\M_v$.

The idea is that a word $w\in u\convW v$ satisfies $\WC(w)=K$ iff
the last $|v|$ letters of $w$ have specific values, depending on $K'$.
Indeed, by definition of $\WC$, if $K'$ is not coarser than $I$,
or if $K''\not=J$, then the coefficient of $U_K$ is zero.
Now, let us fix a composition $I=(i_1,\dots,i_l)$ and a composition $K$
satisfying the conditions of the theorem. Then any word of
$u\convW v=w=u'\cdot v'$ satisfying $\WC(w)=K$ has also $|K|$ different
letters.
Now, for all $j<l$, if $I_j$ and $I_{j+1}$ come from the same part of $K$,
the letter $u'_{i_1+\dots+i_j}$ has to appear in $v'$, otherwise this
letter does not appear in $v'$. Hence, given the letters appearing in $u'$,
the letters appearing in $v'$ are also fixed, which completely determines
$v'$ too (since its packed word is given). The number of ways of choosing the
letters appearing in $u'$ obviously is the binomial coefficient
$\binomial{l(K)}{l(I)}$.
\qed

\begin{example}
{\rm
Let $I=(2,2,1)$, $J=(1,3)$, and $K=(4,1,1,3)$.
We can choose $\sigma=23212$ and $\tau=2111$.
Then $K'=(4,1)$ and $K''=(1,3)$. The coefficient of $K$ in
$U_I U_J$ is $\binom{4}{2}=6$ and, indeed, there are four packed words
in the (modified) convolution $11223\convW 1222$ with W-composition $K$:
\begin{equation}
112241333,\ 113341222,\
112231444,\
223341222.
\end{equation}
Those four packed words are obtained as follows: since $(4,1,1,3)$ is obtained
by gluing together the first two parts of $I$, this means that, if $\WC(w)=K$,
the last four letters of $w$ have to be the first letter of $w$ or the one not
used in its first five letters.
}
\end{example}

\begin{note}
{\rm
This quotient is not a Hopf quotient, since again,
${\cal J}'$ is not a coideal. For example, $221\sim 112$ but
\begin{equation}
\overline\Delta(\M_{221}) = \M_{1}\otimes \M_{11},
\text{\ and\ }
\overline\Delta(\M_{112}) = \M_{11}\otimes \M_{1}.
\end{equation}
}
\end{note}

%\begin{note}
%{\rm
%The proof of Theorem~\ref{TEQuotb} can be transformed in a direct proof of
%Theorem~\ref{TEQuot} since, by definition of the realizations of all those
%algebras,
%\begin{equation}
%\F_\sigma := \sum_{u| \Std(u)=\sigma} M_u.
%\end{equation}
%}
%\end{note}

%%%%%%%%%%%%%%%%%%%%%%%%%%%%%%%%%%%%%%%%%%%%%%%%%%%%%%%%%%%%%%%%%%%%%%%%%%%%%%%
\subsection{Change of bases in $\NCSF$}

As in the case of permutations, we have a map going from $\NCSF$ to itself:
start with the injection of $\NCSF$ into $\WQSym$, which reads
\begin{equation} 
R_I := \sum_{\CDes(u)=I} \M_u,
\end{equation}
and then go from $\WQSym$ to $\NCSF$ by the W-quotient homomorphism.

Let $\phi'$ be the composition of those maps and let $R'_I$ be the image of
$R_I$ by $\phi'$:
\begin{equation}
\begin{array}{cccc}
\phi':& \NCSF & \to     & \TEb \\
     & R_I   & \mapsto & R'_I. \\
\end{array}
\end{equation}
By definition of $\phi'$, we have
\begin{equation}
\label{RpPsi}
R'_I := \sum_{\CDes(u)=I} \overline{\M_u}
=
\sum_{\genfrac{}{}{0pt}{}%
{\CDes(u)=I}{\WC(u)=J}} \Psi_J,
\end{equation}
where  $\overline\M_u$ is the image of $\M_u$ by the W-quotient
homomorphism.
Then, since $\Psi_{1^n}=R_{1^n}$ and
$R'_{1^n}=\overline{\M_{n\dots 21}}=U_{1^n}$, we have
$\phi'(\Psi_n)=U_n$ for all $n$, so that, thanks to the product formulas of
$\Psi_n$ and $U_n$, $\phi'(\Psi_I)=U_I$ for all compositions $I$.

Since the $U_n$ are algebraic generators of $\TEb$, the algebra morphism
$\phi'$ is an isomorphism of algebras, so that, applying $\phi'^{-1}$ to
Equation~(\ref{RpPsi}), one gets

\begin{theorem}
\label{ThRPsi}
Let $I$ be a composition of $n$. Then
\begin{equation}
\label{R2Psi}
R_I = \sum_{J\vDash n} K_{IJ} \Psi_J,
\end{equation}
where $K_{IJ}$ is the number of permutations satisfying
$\CDes(u)=I$ and $\WC(u)=J$.
In particular, the $K_{IJ}$ are nonnegative integers.
\end{theorem}

Combining this last result with Equation~(\ref{LPsi}), one then gets back
Corollary~\ref{CorRPsi}.
Examples of the transition matrices are given in Section~\ref{tabsk},
together with the same matrices filled with the corresponding packed words.

%%%%%%%%%%%%%%%%%%%%%%%%%%%%%%%%%%%%%%%%%%%%%%%%%%%%%%%%%%%%%%%%%%%%%%%%%%%%%%%
%%%%%%%%%%%%%%%%%%%%%%%%%%%%%%%%%%%%%%%%%%%%%%%%%%%%%%%%%%%%%%%%%%%%%%%%%%%%%%%
%%%%%%%%%%%%%%%%%%%%%%%%%%%%%%%%%%%%%%%%%%%%%%%%%%%%%%%%%%%%%%%%%%%%%%%%%%%%%%%
\section{Tables}
\label{tabs}

%%%%%%%%%%%%%%%%%%%%%%%%%%%%%%%%%%%%%%%%%%%%%%%%%%%%%%%%%%%%%%%%%%%%%%%%%%%%%%%
\subsection{Coefficients $G_{IJ}$}
\label{tabsg}

Here are the transition matrices from $R$ to $L$ (the
matrices of the coefficients $G_{IJ}$) for $n=3$ and $n=4$, the compositions
being in lexicographic order.
To save space and for better readability, $0$ has been represented by a dot.

\begin{equation}
M_3(R,L) =
\left(
\begin{matrix}
1 & . & . & . \\
. & 2 & 1 & . \\
. & . & 1 & . \\
. & . & . & 1 
\end{matrix}
\right)
\end{equation}

\begin{equation}
M_4(R,L) =
\left(
\begin{matrix}
1 & . & . & . & . & . & . & . \\
. & 3 & 2 & . & 1 & 1 & . & . \\
. & . & 2 & . & 1 & . & . & . \\
. & . & 1 & 3 & . & 2 & 1 & . \\
. & . & . & . & 1 & . & . & . \\
. & . & . & . & . & 2 & 1 & . \\
. & . & . & . & . & . & 1 & . \\
. & . & . & . & . & . & . & 1
\end{matrix}
\right)
\end{equation}

Here are the same matrices with the list of permutations having a given
recoil composition (or descent composition of the inverse) and G-composition,
instead of their number.
\begin{equation}
{\goth M}_3 =
%\left(
\begin{array}{|c||c|c|c|c|}
\hline
\text{\rm $\GC\backslash$ Rec}& 3 & 21 & 12 & 111 \\[.1cm]
\hline
\hline
3   & 123 &   &   &   \\[.1cm]
\hline
21  &     & \gf{132}{312} & 231 &   \\[.1cm]
\hline
12  &     &   & 213 &   \\[.1cm]
\hline
111 &     &   &   & 321 \\[.1cm]
\hline
\end{array}
%\right)
\end{equation}

\begin{equation}
{\goth M}_4 =
%\left(
\begin{array}{|c||c|c|c|c|c|c|c|c|}
\hline
\text{\rm $\GC\backslash$ Rec}& 4 & 31 & 22 & 211 & 13 & 121 & 112 & 1111
\\[.1cm]
\hline
\hline
4& 1234 &   &   &   &   &   &   &   \\[.1cm]
\hline
31&  & \gf{1243,\ 1423}{4123} & \gf{1342}{3412} &   & 2341 & 2413 &   &  
\\[.1cm]
\hline
22&     &   & \gf{1324}{3124} &   & 2314 &   &   &   \\[.1cm]
\hline
211&     &   & 3142 & \gf{1432,\ 4132}{4312} &   & \gf{2431}{4231} & 3241 &  
\\[.1cm]
\hline
13&     &   &   &   & 2134 &   &   &   \\[.1cm]
\hline
121&     &   &   &   &   & \gf{2143}{4213} & 3421 &   \\[.1cm]
\hline
112&     &   &   &   &   &   & 3214 &   \\[.1cm]
\hline
1111&     &   &   &   &   &   &   & 4321 \\[.1cm]
\hline
\end{array}
%\right)
\end{equation}

%%%%%%%%%%%%%%%%%%%%%%%%%%%%%%%%%%%%%%%%%%%%%%%%%%%%%%%%%%%%%%%%%%%%%%%%%%%%%%%
\subsection{Coefficients $K_{IJ}$}
\label{tabsk}

Here are the transition matrices from $R$ to $\Psi$ (the
matrices of the coefficients $K_{IJ}$) for $n=3$ and $n=4$, the compositions
being in lexicographic order.
To save space and for better readability, $0$ has been represented by a dot.

\begin{equation}
M_3(R,\Psi) =
\left(
\begin{matrix}
1 & . & . & . \\
1 & 2 & 1 & . \\
1 & . & 1 & . \\
1 & 2 & 2 & 1 
\end{matrix}
\right)
\end{equation}

\begin{equation}
M_4(R,\Psi) =
\left(
\begin{matrix}
1 & . & . & . & . & . & . & . \\
1 & 3 & 2 & . & 1 & 1 & . & . \\
1 & . & 2 & . & 1 & . & . & . \\
1 & 3 & 5 & 3 & 2 & 3 & 1 & . \\
1 & . & . & . & 1 & . & . & . \\
1 & 3 & 2 & . & 2 & 3 & 1 & . \\
1 & . & 2 & . & 2 & . & 1 & . \\
1 & 3 & 5 & 3 & 3 & 5 & 3 & 1
\end{matrix}
\right)
\end{equation}

Here are the same matrices with the list of packed words having a given
descent composition and W-composition, instead of their number.
\begin{equation}
{\goth M}'_3 =
%\left(
\begin{array}{|c||c|c|c|c|}
\hline
\text{\rm $\WC\backslash\CDes$}& 3 & 21 & 12 & 111 \\[.1cm]
\hline
\hline
3   & 111 &   &   &   \\[.1cm]
\hline
21  & 112 & \gf{121}{221} & 212 &   \\[.1cm]
\hline
12  & 122 &               & 211 &   \\[.1cm]
\hline
111 & 123 & \gf{132}{231} & \gf{312}{213} & 321 \\[.1cm]
\hline
\end{array}
%\right)
\end{equation}

\begin{equation*}
{\goth M}'_4 =
\begin{array}{|c||c|c|c|c|c|c|c|c|}
\hline
\text{\rm $\WC\backslash\CDes$}& 4 & 31 & 22 & 211 & 13 & 121 & 112 & 1111
\\[.1cm]
\hline
\hline
4   & 1111 &   &   &   &   &   &   &   \\[.1cm]
\hline
31  & 1112 & \gf{1121,\ 1221}{2221} & \gf{2212}{1212} &   & 2112 & 2121 &   &  
\\[.1cm]
\hline
22  & 1122 &   & \gf{1211}{2211} &   & 2122 &   &   &   \\[.1cm]
\hline
211 & 1123 & \gf{1132,1231}{2231} & \gf{1213,1312,2213}{2312,3312}
& \gf{1321,2321}{3321} & \gf{2123}{3123} & \gf{2132,3132}{3231} & 3213 &  
\\[.1cm]
\hline
13  & 1222 &   &   &   & 2111 &   &   &   \\[.1cm]
\hline
121 & 1223 & \gf{1232,1332}{2331} & \gf{1323}{2313} &   & \gf{2113}{3112}
 & \gf{2131,3121}{3221} & 3212 &   \\[.1cm]
\hline
112 & 1233 &   & \gf{1322}{2311} &   & \gf{2133}{3122} &   & 3211 &   \\[.1cm]
\hline
1111& 1234 & \gf{1243,1342}{2341} & \gf{1324,1423,2314}{2413,3412} &
\gf{1432,2431}{3421} & \gf{2134,3124}{4123} & \gf{2143,3142,3241}{4132,4231}
 & \gf{3214,4213}{4312} & 4321 \\
\hline
\end{array}
\end{equation*}

%%%%%%%%%%%%%%%%%%%%%%%%%%%%%%%%%%%%%%%%%%%%%%%%%%%%%%%%%%%%%%%%%%%%%%%%%%%%%%%
%%%%%%%%%%%%%%%%%%%%%%%%%%%%%%%%%%%%%%%%%%%%%%%%%%%%%%%%%%%%%%%%%%%%%%%%%%%%%%%
%%%%%%%%%%%%%%%%%%%%%%%%%%%%%%%%%%%%%%%%%%%%%%%%%%%%%%%%%%%%%%%%%%%%%%%%%%%%%%%

\end{document}